\documentclass{amsart}

\usepackage{amsthm,amsmath,amssymb}
\usepackage[margin=1in]{geometry}
\usepackage{datetime}
\usepackage{graphicx}
\usepackage[numbers]{natbib}

\makeatletter
\patchcmd{\@settitle}{\uppercasenonmath\@title}{}{}{}
\patchcmd{\@setauthors}{\MakeUppercase}{}{}{}
\patchcmd{\section}{\scshape}{}{}{}
\makeatother

\begin{document}

\title{Four monoclinic Diophantine parallelepiped parameterizations}
\author{Randall L. Rathbun}
\email{randallrathbun@gmail.com}

\begin{abstract}
Four parametrizations for the bi-orthogonal monoclinic Diophantine parallelepiped are given.
\end{abstract}

\maketitle

\setcounter{section}{0}
\section*{{\bf The bi-orthogonal monoclinic Diophantine parallelepiped}}

The bi-orthogonal monoclinic Diophantine parellelepiped is a cuboid with 2 right angles at a vertice and
the other angle not a right angle. This figure is composed of two congruent parallelograms
joined by 4 orthogonal rectangles, all of whose edges are rational or integer.
\begin{figure}[!ht]
\centering
\includegraphics[scale=1.0]{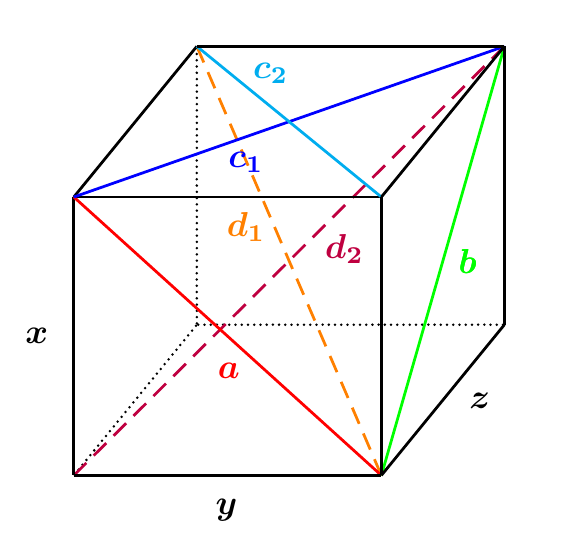}
\vspace{-6pt}
\caption{The Bi-orthogonal Monoclinic Cuboid.}
\label{fig:mono}
\end{figure}
Let $x$, $y$, $z$, denote the three different edges, with $x,y,z \in Z$.
The face rectangle $(x,y)$ has diagonal $a$, the face rectangle $(x,z)$ has diagonal $b$
and the face parallelogram $(y,z)$ has diagonals $c_1, c_2$.
The two different body diagonals are denoted by $d_1, d_2$.

These lengths satisfy the equations
\begin{align}
x^2 + y^2 & = a^2 \\
x^2 + z^2 & = b^2 \\
x^2 + c_1^2 & = d_1^2 \\
x^2 + c_2^2 & = d_2^2 \\
2y^2 + 2z^2 & = c_1^2 + c_2^2 \\
2y^2 + 2b^2 & = d_1^2 + d_2^2 \\
2a^2 + 2z^2 & = d_1^2 + d_2^2
\end{align}

The bi-orthogonal monoclinic piped is a perfect parallelepiped, or a piped which has all twelve edges rational,
all six face diagonals rational, and all four body diagonals rational.

Perfect parallelepipeds have been considered since the question was first raised by Richard Guy in his
UPINT\cite{guy}, Section D. Diophantine equations, Problem D18, {\it "Is there a perfect cuboid?"} where Guy
asks if a perfect parallelepiped can be found, even if no perfect cuboid is found.

Jorge Sawyer and Clifford Reiter discovered a perfect parallelepiped\cite{saw}. Clifford Reiter and Jordan Tirrell
then released an additional series of papers on parameterizing them\cite{reit1}\cite{reit2}. A more recent paper on an infinite
family of perfect parallelepipeds has also appeared\cite{soko}.

Walter Wyss was able to successfully parameterize the perfect parallelogram which has rational sides and two
rational diagonals\cite{wyss1}\cite{wyss2}, and he extended this to parameterize a family of bi-orthogonal monoclinic
pipeds\cite{wyss3}.

Ruslan Sharipov recast Wyss' parametrizations, and extended them to three additional mappings for the bi-orthogonal monoclinic
piped\cite{shar}.

By taking advantage of the polynomial equation (3.9), which governs monoclinic pipeds, in Sharipov's paper\cite{shar},
some new parametrizations of the monoclinic piped have been discovered.

In this paper, four parametrizations of the monoclinic piped are presented where $m,n \in \mathbb{Z}$ are integers.

\section*{{\bf Parametrization 1}}

The first parametrization, $m,n\in\mathbb{Z}$, is:

\begin{align}
\label{eq:parm1}
\begin{split}
x & = | \; 4(n^2-2m^2)(n^2+2mn+2m^2)(n^2+4mn+2m^2)  \; | \\
y & = | \; 3(n^2-2m^2)(n^2+2mn+2m^2)(n^2+4mn+2m^2)  \; | \\
z & = | \; 16mn(n+m)(n+2m)(n^2+2mn+2m^2)  \; | \\
c_1 & = | \; (n^2+4mn+2m^2)(n^2+4mn+6m^2)(3n^2+4mn+2m^2)  \; | \\
c_2 & = | \; (n^2-2m^2)(n^2+2m^2)(3n^2+8mn+6m^2)  \; | \\
d_1 & = | \; (n^2+4mn+2m^2)(5n^4+16mn^3+28m^2n^2+32m^3n+20m^4)  \; | \\
d_2 & = | \; (n^2-2m^2)(5n^4+24mn^3+52m^2n^2+48m^3n+20m^4)  \; | \\
a & = | \; 5(n^2-2m^2)(n^2+2mn+2m^2)(n^2+4mn+2m^2)  \; | \\
b & = | \; 4(n^2+2mn+2m^2)^3  \; |
\end{split}
\end{align}

$\left(\frac{m}{n}\right)$ has to satisfy certain ranges to make a piped: \par
For the positive range: \par \vspace{-0.5em}
\hspace{10mm} let  $r_1 = \sqrt[\leftroot{0}\uproot{1}4]{12x^4 + 24x^3 + 16x^2 + 4x - 3} \approxeq 0.28126795021\dots$ \par \vspace{-0.6em}
\hspace{10mm} let $r_2 = \frac{1}{2r_1} \approxeq 1.77766432195\dots$ \par \vspace{-0.6em}
then \par \vspace{-1em}
\hspace{30mm} $ 0 < \frac{m}{n} < r_1 \quad\text{or}\quad r_2 < \frac{m}{n} < +\infty $ \par \vspace{-0.6em}
for the negative range: \par \vspace{-0.5em}
\hspace{10mm} let $r_3 = -\sqrt[\leftroot{0}\uproot{1}4]{12x^4 + 56x^3 + 80x^2 + 52x + 13} \approxeq -0.60976156477\dots$ \par \vspace{-0.6em}
\hspace{10mm} let $r_4 = -\sqrt[\leftroot{0}\uproot{1}4]{52x^4 + 104x^3 + 80x^2 + 28x + 3} \approxeq -0.81999264776\dots$ \par \vspace{-0.6em}
then \par \vspace{-1em}
\hspace{30mm} $-\frac{1}{2} > \frac{m}{n} > r_3 \quad\text{or}\quad r_4 > \frac{m}{n} > -1$

\section*{{\bf Parametrization 2}}

The second parametrization,, $m,n\in\mathbb{Z}$, is:

\begin{align}
\label{eq:parm2}
\begin{split}
x & = | \; 4(n^2-4mn-4m^2)(5n^2+8mn+4m^2)(7n^2+12mn+4m^2)  \; | \\
y & = | \; 3(n^2-4mn-4m^2)(5n^2+8mn+4m^2)(7n^2+12mn+4m^2)  \; | \\
z & = | \; 32n(n+m)(n+2m)(3n+2m)(5n^2+8mn+4m^2)  \; | \\
c_1 & = | \; (7n^2+12mn+4m^2)(9n^2+20mn+12m^2)(11n^2+12mn+4m^2)  \; | \\
c_2 & = | \; (n^2-4mn-4m^2)(3n^2+4mn+4m^2)(17n^2+28mn+12m^2)  \; | \\
d_1 & = | \; (7n^2+12mn+4m^2)(101n^4+312mn^3+424m^2n^2+288m^3n+80m^4)  \; | \\
d_2 & = | \; (n^2-4mn-4m^2)(149n^4+488mn^3+616m^2n^2+352m^3n+80m^4)  \; | \\
a & = | \; 5(n^2-4mn-4m^2)(5n^2+8mn+4m^2)(7n^2+12mn+4m^2)  \; | \\
b & = | \; 4(5n^2+8mn+4m^2)^3  \; |
\end{split}
\end{align}

$\left(\frac{m}{n}\right)$ has to satisfy certain ranges to make a piped: \par
For the positive range \par \vspace{-0.5em}
\hspace{10mm} let $r_1 = \sqrt[\leftroot{0}\uproot{1}4]{48x^4 + 64x^3 - 40x^2 - 112x - 53} \approxeq 1.27766432195\dots$ \par \vspace{-0.6em}
then \par \vspace{-1em}
\hspace{30mm} $ r_1 < \frac{m}{n} < +\infty $ \par \vspace{-0.6em}
for the negative range \par \vspace{-0.5em}
\hspace{10mm} let $r_2 = -\sqrt[\leftroot{0}\uproot{1}4]{48x^4 + 192x^3 + 280x^2 + 176x + 27} \approxeq -0.21873204978\dots$ \par \vspace{-0.6em}
\hspace{10mm} let $r_3 = -\sqrt[\leftroot{0}\uproot{1}4]{208x^4 + 832x^3 + 1256x^2 + 848x + 213} \approxeq -1.31999264776\dots$ \par \vspace{-0.6em}
then \par \vspace{-1em}
\hspace{30mm} $ r_2 > \frac{m}{n} > -\frac{1}{2} \quad\text{or}\quad r_3 > \frac{m}{n} > -\frac{3}{2}$

\section*{{\bf Parametrization 3}}

The third parametrization, $m,n\in\mathbb{Z}$, is:

\begin{align}
\label{eq:parm3}
\begin{split}
x & = | \; 4(7n^2-8mn-16m^2)(13n^2+24mn+16m^2)(17n^2+40mn+16m^2)  \; | \\
y & = | \; 3(7n^2-8mn-16m^2)(13n^2+24mn+16m^2)(17n^2+40mn+16m^2)  \; | \\
z & = | \; 32n(n+4m)(3n+4m)(5n+4m)(13n^2+24mn+16m^2)  \; | \\
c_1 & = | \; (17n^2+40mn+16m^2)(19n^2+56mn+48m^2)(33n^2+40mn+16m^2)  \; | \\
c_2 & = | \; (7n^2-8mn-16m^2)(9n^2+8mn+16m^2)(43n^2+88mn+48m^2)  \; | \\
d_1 & = | \; (17n^2+40mn+16m^2)(725n^4+2384mn^3+3808m^2n^2+3328m^3n+1280m^4)  \; | \\
d_2 & = | \; (7n^2-8mn-16m^2)(965n^4+3856mn^3+6112m^2n^2+4352m^3n+1280m^4)  \; | \\
a & = | \; 5(7n^2-8mn-16m^2)(13n^2+24mn+16m^2)(17n^2+40mn+16m^2)  \; | \\
b & = | \; 4(13n^2+24mn+16m^2)^3  \; |
\end{split}
\end{align}

$\left(\frac{m}{n}\right)$ has to satisfy certain ranges to make a piped: \par
For the positive range \par \vspace{-0.5em}
\hspace{10mm} let $r_1 = \sqrt[\leftroot{0}\uproot{1}4]{768x^4 + 2304x^3 + 2464x^2 + 1104x - 37} \approxeq 0.0312679502117\dots$ \par \vspace{-0.6em}
\hspace{10mm} let $r_2 = \sqrt[\leftroot{0}\uproot{1}4]{768x^4 + 256x^3 - 1120x^2 - 1328x - 453} \approxeq 1.52766432195\dots$ \par \vspace{-0.6em}
then \par \vspace{-1em}
\hspace{30mm} $ 0 < \frac{m}{n} < r_1 \quad\text{or}\quad r_2 < \frac{m}{n} < +\infty$ \par \vspace{-0.6em}
for the negative range \par \vspace{-0.5em}
\hspace{30mm} $ 0 > \frac{m}{n} > -\frac{1}{4} \quad\text{or}\quad -\frac{3}{4} > \frac{m}{n} > -\frac{5}{4}$

\section*{{\bf Parametrization 4}}

The fourth parametrization, $m,n\in\mathbb{Z}$, is:

\begin{align}
\label{eq:parm4}
\begin{split}
x & = | \; 4(7n^2-12mn-18m^2)(17n^2+30mn+18m^2)(23n^2+48mn+18m^2)  \; | \\
y & = | \; 3(7n^2-12mn-18m^2)(17n^2+30mn+18m^2)(23n^2+48mn+18m^2)  \; | \\
z & = | \; 48n(n+3m)(4n+3m)(5n+6m)(17n^2+30mn+18m^2)  \; | \\
c_1 & = | \; 9(3n^2+8mn+6m^2)(23n^2+48mn+18m^2)(41n^2+48mn+18m^2)  \; | \\
c_2 & = | \; 3(7n^2-12mn-18m^2)(11n^2+12mn+18m^2)(19n^2+36mn+18m^2)  \; | \\
d_1 & = | \; (23n^2+48mn+18m^2)(1205n^4+3912mn^3+5940m^2n^2+4752m^3n+1620m^4)  \; | \\
d_2 & = | \; (7n^2-12mn-18m^2)(1685n^4+6288mn^3+9180m^2n^2+6048m^3n+1620m^4)  \; | \\
a & = | \; 5(7n^2-12mn-18m^2)(17n^2+30mn+18m^2)(23n^2+48mn+18m^2)  \; | \\
b & = | \; 4(17n^2+30mn+18m^2)^3  \; |
\end{split}
\end{align}

$\left(\frac{m}{n}\right)$ has to satisfy certain ranges to make a piped: \par
For the positive range \par \vspace{-0.5em}
\hspace{10mm} let $r_1 = \sqrt[\leftroot{0}\uproot{1}4]{324x^4 + 216x^3 - 432x^2 - 636x - 241} \approxeq 1.44433098861\dots$ \par \vspace{-0.6em}
then \par \vspace{-1em}
\hspace{30mm} $r_1 < \frac{m}{n} < +\infty$ \par
for the negative range \par \vspace{-0.5em}
\hspace{10mm} let $ r_2 = -\sqrt[\leftroot{0}\uproot{1}4]{324x^4 - 1080x^3 + 1296x^2 - 660x + 31} \approxeq -0.052065383121\dots$ \par \vspace{-0.6em}
\hspace{10mm} let $ r_3 = -\sqrt[\leftroot{0}\uproot{1}4]{108x^4 + 648x^3 + 1296x^2 + 1132x + 373} \approxeq -0.94309489810\dots$ \par \vspace{-0.6em} 
\hspace{10mm} let $ r_4 = -\sqrt[\leftroot{0}\uproot{1}4]{4212x^4 + 14040x^3 + 17712x^2 + 10020x + 2083} \approxeq -1.15332598109\dots$ \par \vspace{-0.6em}
then \par \vspace{-1em}
\hspace{30mm} $r_2 > \frac{m}{n} > -\frac{1}{3} \quad\text{or}\quad -\frac{5}{6} > \frac{m}{n} > r_3 \quad\text{or}\quad r_4 > \frac{m}{n} > -\frac{4}{3}$

\noindent \rule[0pt]{3.0in}{0.25pt}
\end{document}